\documentclass[usletter]{article}

\usepackage{amsmath, amssymb}
\usepackage{graphicx}

\newtheorem{theorem}{Theorem}
\newtheorem{corollary}[theorem]{Corollary}

\newcommand{\sector}[1]{
  \refstepcounter{section}
  \setcounter{equation}{0}
  \setcounter{theorem}{0}
  \mbox{\ }\smallskip
  {\begin{center}
    \bf \thesection.\ {#1}
  \end{center}}}

\newcommand{\example}{\refstepcounter{theorem}
   \vspace*{\medskipamount}
   {\noindent\bf Example \thesection.\arabic{theorem}}\mbox{\ \ }}  
\newcommand{\eofproof}{{ \hspace{0.2cm}\\
        \mbox{\,}\hfill\framebox[0.3cm]{\rule{0cm}{.1cm} } }}
\newcommand{\Pf}[1]{{\noindent \bf Proof{#1}\ }}

\newcommand{\eqdef}{{\ \stackrel{\mathrm{def}}{=}\ }}

\begin{document}

\title{Asymptotic analysis of a size-structured cannibalism model with infinite dimensional environmental
feedback}
\author{J\'{o}zsef Z. Farkas\thanks{E-mail: jzf@maths.stir.ac.uk (corresponding author)}\ {\normalsize $^1$}\mbox{\ }
and Thomas Hagen\thanks{E-mail: thagen@memphis.edu}\ {\normalsize $^2$} \mbox{\
}
    \mbox{\ }\\
    \mbox{\ }\\
    {\normalsize $^1$\,Department of Computing Science and Mathematics,
     University of Stirling,}\\
    {\normalsize  Stirling, FK9 4LA, UK}  \\
{\normalsize $^2$\,Department of Mathematical Sciences, The University of Memphis,}\\
    {\normalsize  Memphis, TN 38152, USA} }
\date{\normalsize September 2008\\ in revised form: December 2008}

\maketitle

\begin{abstract}
In this work we consider a 
size-structured cannibalism model with the model ingredients
(fertility, growth, and mortality rate) depending on size 
(ranging over an infinite domain) and on a
general function of the standing population (environmental
feedback). Our focus is on the asymptotic behavior of the
system, in particular on the effect of cannibalism on the long-term
dynamics. To this end, we formally linearize the system about
steady state and establish conditions in terms of the model
ingredients which yield uniform exponential stability of the governing linear semigroup. 
We also show how the point spectrum of the linearized semigroup generator
can be characterized in the special case of a separable attack
rate and establish a  general instability result.
Further spectral analysis allows us  to give conditions for 
asynchronous exponential growth of the linear semigroup.

\end{abstract}

{\samepage
\sector{Introduction}

Cannibalism is a phenomenon  observed among many species, e.g.\,certain
fish populations. Sophisticated population models are capable to
 elucidate a potentially stabilizing effect of cannibalism, underscoring 
that certain populations may benefit from cannibalism
when resources are limited. Consequently, the effects of cannibalism on the long-term
dynamics of populations have attracted considerable interest and have 
been analyzed for various structured
population models (see \cite{CUS,GDR} for further references).
}

Structured population models are typically formulated as partial
differential equations for population densities. Diekmann et
al.\,have developed a  general mathematical framework to study
analytical questions for structured populations (see
\cite{D1,D2,DGM2}), including those pertaining to linear/nonlinear
stability of population equilibria. In this context it was
recently proven for large classes of structured population models,
formulated as integral (or delay) equations,
that  the nonlinear stability/instability of a population
equilibrium is completely determined by its linear
stability/instability, a result commonly referred to  as the
``Principle of Linearized Stability''.

In a series of recent papers we have successfully applied linear
semigroup { and spectral} methods to formulate biologically interpretable
conditions for the linear stabili\-ty/in\-sta\-bi\-li\-ty of
equilibria of various structured population models, see
\cite{FH3,FH1} and the references therein. In these problems we assumed that any effect of
intraspecific competition {between individuals of different sizes} 
on individual behavior is primarily due
to a change in population size and that every individual in the
population can influence the vital rates of other individuals
(``scramble competition").

When the competition among individuals is based upon some
hierarchy in the population, often related to the size of
individuals, environmental feedback is incorporated in the model
through infinite dimensional interaction variables (``contest
competition"). This situation in its simplest form is of relevance
for a forest consisting of tree individuals in which the height of
a tree determines its rank in the population. Taller individuals
have higher efficiency when competing for resources such as light,
while individuals of lower rank do not influence the vital rates
of individuals of higher rank. 

Of interest in this work is the linear stability analysis of
population equilibria  of a continuous quasilinear size-structured
model, recently discussed in \cite{D1}. In particular, it was
proven there that the Principle of Linearized Stability holds true
for this model. The density evolution of individuals of size $s$
is assumed to be governed by the partial differential equation
\begin{equation}
n_t(s,t)+\left(\gamma(s,E(s,t))\,n(s,t)\right)_s+\left(\mu(s,E(s,t))+M(s,t)\right)\,n(s,t)=0,\label{pde}
\end{equation}
defined for $s\in (0,\infty)$ and $t>0$. The density of zero
(or minimal) size individuals is given by the nonlocal boundary
condition
\begin{equation}
\gamma(0,E(0,t))\,n(0,t)=\int_0^{\infty}
\beta(s)\,n(s,t)\,ds,\quad t>0. \label{boundary}
\end{equation}
The initial condition reads
\begin{equation}
n(s,0)=n_0(s),\quad   s\in [0,\infty). \label{initial}
\end{equation}
Here $\beta$, $\mu$ and $\gamma$ denote the fertility, mortality
and growth rate of individuals, respectively. All vital rates are
size-dependent. Moreover, it is assumed that the mortality and
growth rates of individuals depend on the function
\begin{equation}
E(s,t)=\int_0^{\infty}c(y)\,\alpha(y,s)\,n(y,t)\,dy.\label{density}
\end{equation}
Here $\alpha$ denotes the size-specific attack rate, while $c$
denotes the energetic value of the attacked prey. Hence the
function $E$ models the assumption that the extra energy
available  due to cannibalism is channeled into the growth of
individuals and affects their starvation-driven mortality. The
extra size-specific mortality rate due to cannibalism is given by
\begin{equation}\label{density2}
M(s,t)=\int_0^{\infty}\alpha(s,y)\,n(y,t)\,dy.
\end{equation}
(Note the switch of variables in $\alpha$ in contrast to
\eqref{density}\,.)
We remark that the cannibalism model discussed here incorporates
 the environmental feedback variable used in the hierarchical size-structured model in \cite{FH4}.

We impose the following regularity conditions on the model ingredients:
\begin{itemize}
\item $\mu=\mu(s,E)\in C_b([0,\infty);C^1_b[0,\infty))$, $\mu\geq 0$
\item $\gamma=\gamma(s,E)\in C^1_b([0,\infty); C^1_b[0,\infty))$, $\gamma\geq
\gamma_0>0$ for some constant $\gamma_0$
\item $\beta=\beta(s)\in C_b([0,\infty))$, $\beta\geq 0$
\item $c=c(s)\in C_b([0,\infty))$, $c\geq 0$
\item $\alpha=\alpha(y,s)\in C_b([0,\infty); C^1_b([0,\infty)))$, $\alpha\geq 0$.
\end{itemize}
The subscript $b$ indicates that functions (and derivatives in
case of $C^1$) are bounded. For clarity in later developments we will write $D_2\alpha$ for the derivative of $\alpha$ with respect to its second argument. The regularity assumptions above are tailored toward
the linear analysis of this work. They might, however, not suffice
to guarantee the existence and uniqueness of solutions of
Eqs.~\eqref{pde}--\eqref{density2}, even in the steady-state case. Well-posedeness of structured
partial differential equation models with infinite dimensional
environmental feedback variables is in general an open question.
It has recently been shown in \cite{AI} that population models
with infinite dimensional interaction variables may exhibit a more
complicated dynamical behavior than the simple size-structured
model of scramble competition.

The study of hierarchical models in the literature (see e.g.\
\cite{CS2} and the references therein) is largely based on a
decoupling of the total population quantity from the governing
equations and a transformation  of the nonlocal partial
differential equation \eqref{pde} into a local one. This technique
allows to prove well-posedness and to study the asymptotic
behavior of solutions by means of ODE methods. For
Eqs.~\eqref{pde}--\eqref{density2} this approach fails, however,
since the mortality and growth rates depend on both size $s$ and
on the environment $E(s,t)$. Therefore it seems unavoidable to
study the original partial differential equation \eqref{pde} with
the nonlocal integral boundary condition \eqref{boundary}
directly.
Moreover, studying the linear stability of stationary solutions for models with infinite
dimensional interaction variables by spectral analysis (as done previously for simpler cases in \cite{FH3,FH1,PR1}) 
has proven  difficult since eigenvalues are not given by an explicitly available characteristic equation (see \cite{FH4}).
Therefore we devised a different approach in \cite{FH4} to { establish 
a linear stability condition} for a particular model with infinite
dimensional interaction variables.

Clearly, system \eqref{pde}-\eqref{density2} admits the trivial solution
$n\equiv 0$. In general, system \eqref{pde}-\eqref{density2}
yields for a stationary solution $n_*$
\begin{equation}\label{statsol}
n_*(s)={{n_*(0)\,\gamma(0,E_*(0))}\over {\gamma(s,E_*(s))}}\,\exp\left\{-\int_0^s\frac{\mu(y,E_*(y))+M_*(y)}{\gamma(y,E_*(y))}\,dy\right\},
\end{equation}
where
\begin{equation}
E_*(s)=\int_0^{\infty}c(y)\,\alpha(y,s)\,n_*(y)\,dy,\quad
M_*(s)=\int_0^{\infty}\alpha(s,y)\,n_*(y)\,dy.
\end{equation}
Here and later on, starred quantities are stationary counterparts
of  the time-dependent functions in
Eqs.~\eqref{pde}--\eqref{density2}. For obvious reasons, we shall
exclusively consider positive stationary solutions of the form
\eqref{statsol} or the trivial solution $n_*\equiv 0$ in the
following. Moreover, to be consistent with later developments, we
shall assume throughout that stationary solutions $n_*$ have the
regularity $W^{1,1}(0,\infty)$. 

For fixed model ingredients $\beta$, $\mu$, $\gamma$,
$\alpha$, $c$, a positive stationary solution $n_*$ with $E_*$ and
$M_*$ as before satisfies the equation
\begin{equation}\label{netgrowth}
    \int_0^{\infty}\frac{\beta(s)}{\gamma(s,E_*(s))}\exp\left\{-\int_0^s\frac{\mu(y,E_*(y))+M_*(y)}{\gamma(y,E_*(y))}\,dy\right\}\,ds=1.
\end{equation}
Therefore, for solutions $n=n(s,t)$ of the governing equations, it is natural to
define the functional
\begin{equation}\label{netrep}
R(n)=\int_0^{\infty}\frac{\beta(s)}{\gamma(s,E(s,\cdot))}\exp\left\{-\int_0^s\frac{\mu(y,E(y,\cdot))+M(y,\cdot)}{\gamma(y,E(y,\cdot))}\,dy\right\}\,ds,
\end{equation}
where $n$, of course, determines the quantities $E$ and $M$ via \eqref{density}, \eqref{density2}.
$R$ may be regarded as the net reproduction rate of the standing population. 
Note that Eq.~\eqref{netgrowth} requires $R(n_*)=1$ for any positive stationary solution $n_*$.

In \cite{DGM2} a general framework was developed to establish the existence of steady states for
general physiologically structured population models, however with finite dimensional interaction variables. 
For models with infinite dimensional interaction variables one can usually not formulate 
elegant necessary and sufficient conditions for the existence of steady state solutions. 
In this situation, one can construct positive stationary solutions by perturbation of
solutions with $\alpha\equiv \text{const}$. In the latter case, the interaction variables $E$ and $M$ are constant, hence the results in \cite{DGM2} apply.
For an alternative approach we refer the reader to \cite{BCC}.
Throughout the rest of the paper we will tacitly assume that 
stationary solutions of the required regularity are available.

\sector{Stability via dissipativity}

Given a stationary
solution $n_*$, we formally linearize the governing equations by
introducing the infinitesimal perturbation $u=u(s,t)$ and  making
the ansatz $n=u+n_*$. After inserting this expression in the
governing equations and omitting all nonlinear terms, we obtain
the linearized problem
\begin{align}
& u_t(s,t)+\left(\gamma^*(s)\,u(s,t)+\gamma_E(s,E_*(s))\,n_*(s)\,F(s,t)\right)_s+\mu(s,E_*(s))\,u(s,t)
\nonumber \\
& \quad+\mu_E(s,E_*(s))\,n_*(s)\,F(s,t)+M_*(s)\,u(s,t)+n_*(s)\,N(s,t)=0,  \label{linear} \\
&
\gamma^*(0)\,u(0,t)=\int_0^{\infty}\big(\beta(s)-\gamma_E(0,E_*(0))\,n_*(0)\,c(s)\,\alpha(s,0)\big)\,u(s,t)\,ds,
\label{linboundary}
\end{align}
where we have set
\begin{equation}
 F(s,t)=\int_0^{\infty}c(y)\,\alpha(y,s)\,u(y,t)\,dy,\quad
 N(s,t)=\int_0^{\infty}\alpha(s,y)\,u(y,t)\,dy,
\end{equation}
and
\begin{equation}
\gamma^*(s)=\gamma(s,E_*(s)),\quad s\in [0,\infty).
\end{equation}

We denote the Lebesgue space $L^1(0,\infty)$ with its usual
norm $\|\cdot\|$ by ${\mathcal X}$ and introduce the bounded linear
functional $\Lambda$ on $\mathcal X$ by
\begin{equation}
\Lambda(u)=
\int_0^{\infty}\left(\frac{\beta(s)}{\gamma^*(0)}-\frac{\gamma_E(0,E_*(0))\,n_*(0)}{\gamma^*(0)}\,c(s)\,
\alpha(s,0)\right)\,u(s)\,ds.\label{boundaryop}
\end{equation}
Next we define the operators
\begin{align}
    {\mathcal A} u& =-\gamma^*(\cdot)\,u_s,\quad
    \text{Dom}({\mathcal A}) =
    \left\{u\in W^{1,1}(0,\infty)\,|\,u(0)=\Lambda(u)\right\},\label{defA}\\
     {\mathcal B} u & = -\left(\mu(\cdot,E_*(.))+\gamma^*_s(\cdot)+M_*(.)\right)\,u=-\rho^*(.)\,u\quad
    \text{on ${\mathcal X}$,}\label{defB}\\
     {\mathcal C} u&=  -\int_0^{\infty} u(y)\,\left[ c(y)\,\alpha(y,.)\,\left((\gamma_E(.,E_*(.))\,n_*(.))_s+\mu_E(\cdot,E_*(.))\,n_*(.)\right)\right.
\nonumber\\
    &\left.\quad +c(y)\,D_2\alpha(y,.)\,\gamma_E(.,E_*(.))\,n_*(.)+\alpha(.,y)\,n_*(.)\right]\,dy \quad \text{on $\mathcal{X}$} \label{defC}.
\end{align}
Our regularity assumptions on the model ingredients and the
stationary solution guarantee that these operators are
well-defined, that the operators ${\mathcal B}$ and ${\mathcal C}$ are bounded 
on $\mathcal X$, and that the operator ${\mathcal A}$ is closed and densely defined on $\mathcal X$. Thus, the
linearized system \eqref{linear}, \eqref{linboundary} assumes the
form of an initial value problem for an ordinary differential
equation on ${\mathcal X}$
\begin{equation}\label{abstr}
    {d\over dt}\, u = \left({\mathcal A} + {\mathcal B}+ {\mathcal C}\right)\,u,
\end{equation}
with the initial condition
\begin{equation}\label{abstr2}
    u(0)=u_0.
\end{equation}

Our objective in the following is to apply the Lumer-Phillips Theorem from linear semigroup theory (see \cite{Pa}). To obtain stability by virtue of this result, 
we will extend our approach in \cite{FH4}, which was previously
devised for elongational flows in \cite{Ha2}.

\begin{theorem}\label{dissipative}
The operator $\mathcal{A+B+C}$ is the infinitesimal generator of a quasi-contraction semigroup
$T=\{{\mathcal T}(t)\}_{t\geq 0}$ of bounded linear operators on $\mathcal X$. The semigroup is
uniformly exponentially stable if 
\begin{align}
    &\mu(s,E_*(s))+M_*(s)>\big\vert\beta(s)-\gamma_E(0,E_*(0))\,n_*(0)\,c(s)\,\alpha(s,0)
\big\vert \nonumber\\
    &\quad +\int_0^{\infty} \big\vert c(s)\,\alpha(s,y)\,\left[(\gamma_E(y,E_*(y))\,n_*(y))_y+\mu_E(y,E_*(y))\,n_*(y)\right]\nonumber\\
    &\quad +c(s)\,D_2\alpha(s,y)\,\gamma_E(y,E_*(y))\,n_*(y)+\alpha(y,s)\,n_*(y)\big\vert\,dy,\quad
s\geq 0.\label{disscond}
\end{align}
\end{theorem}

\Pf{.} Assume that, for given $h\in \mathcal{X}$ and fixed
$\kappa\in {\mathbb R}$, $u\in \text{Dom}({\mathcal A})$ is such
that, for some $\lambda>0$,
\begin{equation}
    u-\lambda\,({\mathcal A} +{\mathcal B}+ {\mathcal C}+\kappa\,{\mathcal I})\,u=h.\label{res0}
\end{equation}
Then we have
\begin{align}
    \|u\| =&\ \int_0^{\infty} u(s)\,\text{sgn}\,u(s)\,ds\nonumber\\
    =&\ \int_0^{\infty} h(s)\,\text{sgn}\,u(s)\,ds-\lambda\,\int_0^{\infty} (\gamma^*(s)\,u(s))_s\,\text{sgn}\,u(s)\,ds\nonumber\\
    &\ +\lambda\,\int_0^{\infty} (\kappa-\mu(s,E_*(s))-M_*(s))\,u(s)\,\text{sgn}\,u(s)\,ds\nonumber\\
    & -\lambda\,\int_0^{\infty}\int_0^{\infty}u(y)[c(y)\,\alpha(y,s)\,((\gamma_E(s,E_*(s))\,n_*(s))_s+\mu_E(s,E_*(s))\,n_*(s))\nonumber\\
    & \quad +c(y)\,D_2\alpha(y,s)\,\gamma_E(s,E_*(s))\,n_*(s)+\alpha(s,y)\,n_*(s)]\,dy\,\, \text{sgn}\,u(s)\,ds.\label{norm}
\end{align}
First by changing the order of integration, we obtain
\begin{align}
& -\lambda\,\int_0^{\infty}\int_0^{\infty}u(y)\,[c(y)\,\alpha(y,s)\,((\gamma_E(s,E_*(s))\,n_*(s))_s+\mu_E(s,E_*(s))\,n_*(s))\nonumber\\
    & \quad +c(y)\,D_2\alpha(y,s)\,\gamma_E(s,E_*(s))\,n_*(s)+\alpha(s,y)\,n_*(s)]\,dy\,\, \text{sgn}\,u(s)\,ds\nonumber\\
 \le & \ \lambda\,\int_0^{\infty} \vert u(s)\vert \int_0^{\infty}[c(s)\,\alpha(s,y)\,((\gamma_E(y,E_*(y))\,n_*(y))_y+\mu_E(y,E_*(y))\,n_*(y))\nonumber\\
    & \quad +c(s)\,D_2\alpha(s,y)\,\gamma_E(y,E_*(y))\,n_*(y)+\alpha(y,s)\,n_*(y)]\,dy\,ds. \label{ineq1}
\end{align}
Then we note that
the set of points in the interval $(0,\infty)$ where $u$ is
nonzero is the countable union of disjoint open intervals
$(a_i,b_i)$ with $a_i\in [0,\infty)$ and $b_i\in (0,\infty]$ such
that  on each of these intervals either $u>0$ or $u<0$ holds true
with $u(a_i)=0$ unless $a_i=0$, and  $u(b_i)=0$ unless
$b_i=\infty$. If $(a_i,b_i)$ is any such interval on which $u>0$
we have
\begin{align}
& \int_{a_i}^{b_i} h(s)\,\text{sgn}\,u(s)\,ds-\lambda\,\int_{a_i}^{b_i}(\gamma^*(s)\,u(s))_s\,\text{sgn}\,u(s)\,ds\nonumber\\
&\ +\lambda\,\int_{a_i}^{b_i} (\kappa-\mu(s,E_*(s))-M_*(s))\,u(s)\,\text{sgn}\,u(s)ds\nonumber\\
\le & \int_{a_i}^{b_i}\vert h(s)\vert\,ds+\lambda\int_{a_i}^{b_i}\left(\kappa-\mu(s,E_*(s))-M_*(s)\right)u(s)\,ds\nonumber\\
& +\lambda\,\gamma^*(a_i)\,u(a_i)\label{ineq2}.
\end{align}
Similarly, on any interval $(a_i,b_i)$ where $u<0$ we have
\begin{align}
& \int_{a_i}^{b_i} h(s)\,\text{sgn}\,u(s)\,ds-\lambda\,\int_{a_i}^{b_i}(\gamma^*(s)\,u(s))_s\,\text{sgn}\,u(s)\,ds\nonumber\\
&\ +\lambda\,\int_{a_i}^{b_i} (\kappa-\mu(s,E_*(s))-M_*(s))\,u(s)\,\text{sgn}\,u(s)ds\nonumber\\
\le & \int_{a_i}^{b_i}\vert h(s)\vert\,ds+\lambda\int_{a_i}^{b_i}\left(\kappa-\mu(s,E_*(s))-M_*(s)\right)\,\vert u(s)\vert\,ds\nonumber\\
& -\lambda\,\gamma^*(a_i)\,u(a_i)\label{ineq3}.
\end{align}
Combining inequalities \eqref{ineq1}-\eqref{ineq3} we get the estimate
\begin{align}
& \|u\|=  \sum_{i}\int_{a_i}^{b_i} u(s)\,\text{sgn}\,u(s)\,ds\le \|h\|+\lambda\,\gamma^*(0)\,\vert u(0)\vert\nonumber \\
& +\lambda\,\int_0^{\infty} \vert u(s)\vert\ \Big(\kappa-\mu(s,E_*(s))-M_*(s)+\int_0^{\infty}\big\vert c(s)\,\alpha(s,y)\nonumber\\
& \times\big[(\gamma_E(y,E_*(y))\,n_*(y))_y+\mu_E(y,E_*(y))\,n_*(y)\big]\nonumber\\
    & \quad +c(s)\,D_2\alpha(s,y)\,\gamma_E(y,E_*(y))\,n_*(y)+\alpha(y,s)\,n_*(y)\big\vert\,dy \Big)\,ds.
\end{align}
Next we note that
\begin{equation}\label{ineq4}
\vert
u(0)\vert=\vert\Lambda(u)\vert\le\int_0^{\infty}\left\vert\frac{\beta(s)}{\gamma^*(0)}-\frac{\gamma_E(0,E_*(0)n_*(0)}{\gamma^*(0)}c(s)\alpha(s,0)\right\vert
\vert u(s)\vert\,ds.
\end{equation}
Now choose $\kappa\in {\mathbb R}$ such that, for $s\in[0,\infty)$, 
\begin{align}
\kappa \leq & \  \mu(s,E_*(s))+M_*(s)-\left\vert\beta(s)-\gamma_E(0,E_*(0))\,n_*(0)\,c(s)\,\alpha(s,0)
\right\vert\nonumber\\
& \ -\int_0^{\infty}\big\vert c(s)\,\alpha(s,y)\,\big[(\gamma_E(y,E_*(y))\,n_*(y))_y+\mu_E(y,E_*(y))\,n_*(y)\big]\nonumber\\
    & \quad +c(s)\,D_2\alpha(s,y)\,\gamma_E(y,E_*(y))\,n_*(y)+\alpha(y,s)\,n_*(y)\big\vert\,dy.
\end{align}
For such $\kappa$, we have the desired inequality
\begin{equation}
    \|u\|\leq \|h\|\quad\text{for $\lambda>0$,}
\end{equation}
thus establishing dissipativity.
We observe that the operator $\mathcal{A+B+C+\kappa\, I}$ is
densely defined and that the equation
\begin{equation}
    (\lambda I-\mathcal{A})\,u=f
\end{equation}
for $f\in\mathcal{X}$ and $\lambda>0$ sufficiently large has a
unique solution $u\in \text{Dom}\,{\mathcal A}$, given by
\begin{align}
u(s)= & \exp\left\{-\int_0^s\frac{\lambda}{\gamma^*(y)}\,dy\right\}\nonumber\\
&
\times\left(\Lambda(u)+\int_0^s\exp\left\{\int_0^r\frac{\lambda}{\gamma^*(y)}\,dy\right\}\,\frac{f(r)}{\gamma^*(r)}\,dr\right)\label{udef1}
\end{align}
with
\begin{align}
    \Lambda(u) =&
    \left(1-\Lambda\left(\exp\left\{-\int_0^\cdot\frac{\lambda}{\gamma^*(y)}\,dy\right\}\right)\right)^{-1}\,\nonumber\\
    &\Lambda\left(\int_0^\cdot\exp\left\{\int_0^r\frac{\lambda}{\gamma^*(y)}\,dy-\int_0^\cdot
    \frac{\lambda}{\gamma^*(y)}\,dy\right\}\,\frac{f(r)}{\gamma^*(r)}\,dr\right).\label{udef2}
\end{align}
The fact that $u\in \text{Dom}\,{\mathcal A}$ is well defined by
\eqref{udef1}, \eqref{udef2}  follows immediately from the
regularity of the functions involved and their growth behavior. Since $\mathcal B+C+\kappa\,I$
is bounded, the Lumer-Phillips Theorem gives that the operator
$\mathcal{A+B+C+\kappa\, I}$ generates a quasi-contraction semigroup. Specifically, the semigroup $T=\{\mathcal{T}(t)\}_{t\ge 0}$,
generated by the operator $\mathcal{A+B+C}$, satisfies
\begin{equation}
\|\mathcal{T}(t)\|\le e^{-\kappa t},\quad t\geq 0.
\end{equation}
Finally, if condition \eqref{disscond} holds, we can choose $\kappa>0$. Hence the semigroup $T=\{\mathcal{T}(t)\}_{t\ge 0}$ is uniformly exponentially stable.
\eofproof

\begin{corollary}
	A stationary solution $n_*$ of Eqs.~\eqref{pde}--\eqref{density2}
is linearly asymptotically stable if condition \eqref{disscond} holds true.
\end{corollary}

\begin{remark}
For the stability of the trivial equilibrium $n_*\equiv 0$, the
criterion \eqref{disscond} reduces to
\begin{equation}
\mu(s,0)>\beta(s),\quad s\in [0,\infty).\label{trivial}
\end{equation}
Since
\begin{align}
R(0) & =\int_0^{\infty}\frac{\beta(s)}{\gamma(s,0)}\exp\left\{-\int_0^s\frac{\mu(y,0)+0}{\gamma(y,0)}dy\right\}\,ds\nonumber \\
& <
\int_0^{\infty}\bar{\mu}(s,0)\,\exp\left\{-\int_0^s\bar{\mu}(y,0)\,dy\right\}\,ds=1
\end{align}
with $\bar{\mu}(s,0)=\mu(s,0)/\gamma(s,0)$, \eqref{trivial} implies $R(0)<1$. This is the
well-known stability criterion of the trivial steady state in
scramble competition. In the case of the hierarchical model
discussed in \cite{FH4} our stability condition for the trivial
steady state implied $R(0)<1$ as well. 
\end{remark}

\sector{Instability via eigenvalues}

In this section we will characterize part of the point spectrum of the linear semigroup generator when the attack rate $\alpha$ assumes a special form. This result allows us to  establish an instability result, thus complementing our stability result in the previous section.

Throughout this section we make the assumption that the attack rate is separable, i.e.
\begin{equation}
\alpha(x_1,x_2)=\alpha_1(x_1)\,\alpha_2(x_2), \qquad
(x_1,x_2)\in [0,\infty)\times[0,\infty).\label{alpha}
\end{equation}
We can interpret $\alpha_2(x_2)$ as a measure for the likelihood that individuals of size $x_2$ attack, while 
$\alpha_1(x_1)$ represents the likelihood of being attacked at size $x_1$. 
Note that condition \eqref{alpha} assumes no correlation between these two events.  
Even though this assumption might appear as unsatisfactory from a biological point of view, it makes, however, analytical progress possible and henceforth has the potential to shed light on the case of a general attack rate. The main difficulty in case of general attack rates is that the operator $\mathcal{C}$ need not be a finite rank operator. However, our characterization  of the point spectrum of the linearized operator $\mathcal{A+B+C}$ in the following developments relies essentially on this fact. 

We also make the biologically plausible assumption that there is a constant $\mu_0>0$ such that
\begin{equation}\label{muassume}
	\mu(s,E_*(s))+M_*(s)\geq \mu_0\quad\text{for all $s\in[0,\infty)$.}
\end{equation}
Clearly, this assumption is automatically guaranteed if we suppose that $\mu\geq \mu_0$.

The particular choice of the attack rate \eqref{alpha} allows us to cast the operator $\mathcal{C}$ in the form
\begin{equation}\label{newC}
\mathcal{C} u=-\bar{u}_1\,g_1-\bar{u}_2\,g_2,
\end{equation}
where we define
\begin{align}
& \bar{u}_1=\int_0^{\infty}c(s)\,\alpha_1(s)\,u(s)\,ds,\quad \bar{u}_2=\int_0^{\infty}\alpha_2(s)\,u(s)\,ds,\nonumber\\
& g_1=\alpha_2\,\big((\gamma_E\,n_*)_s+\mu_E(.,E_*)\,n_*\big)+\alpha'_2\,\gamma_E(.,E_*)\,n_*,\nonumber\\
& g_2=\alpha_1\,n_*.\label{variables1}
\end{align}
Hence the eigenvalue problem 
\begin{equation*}
(\mathcal{A+B+C}-\lambda\,\mathcal{I})\,u=0,\quad u(0)=\Lambda(u)
\end{equation*}
assumes the form
\begin{align}
& \gamma^*(s)\,u_s+(\rho^*(s)+\lambda)\,u(s)+g_1(s)\,\bar{u}_1+g_2(s)\,\bar{u}_2=0,\label{eigv}\\
& u(0)=\int_0^{\infty}\frac{\beta(s)}{\gamma^*(0)}\,u(s)\,ds-g_3\,\bar{u}_1,\label{eigvbound}
\end{align}
with
\begin{equation*}
g_3=\gamma_E(0,E_*(0))\,n_*(0)\,\alpha_2(0).
\end{equation*}
For 
\begin{equation}\label{lambdacond}
	\text{Re}\,\lambda>-\mu_0,
\end{equation}
the solution of equation \eqref{eigv} is found to be
\begin{align}
u(s)=&\ u(0)\,\pi(s,\lambda)-\bar{u}_1\,\pi(s,\lambda)\,\int_0^s\frac{g_1(r)}{\gamma^*(r)\,\pi(r,\lambda)}\,dr\nonumber\\
&\ -\bar{u}_2\,\pi(s,\lambda)\,\int_0^s\frac{g_2(r)}{\gamma^*(r)\,\pi(r,\lambda)}\,dr,
\label{eigvsol}
\end{align}
where we have made use of the notation
\begin{align}
\pi(s,\lambda)&= \exp\left\{-\int_0^s\frac{\rho^*(y)+\lambda}{\gamma^*(y)}\,dy\right\}\nonumber\\
&=
{{\gamma^*(0)}\over {\gamma^*(s)}}\,\exp\left\{-\int_0^s {{\mu(y,E_*(y))+M_*(y)+\lambda}\over {\gamma^*(y)}}\,dy\right\}.\label{pidef}
\end{align}
Note that condition \eqref{lambdacond} ensures that $u$, given by Eq.~\eqref{eigvsol}, is in $W^{1,1}(0,\infty)$.
We multiply equation \eqref{eigvsol} by $c(s)\,\alpha_1(s)$ and $\alpha_2(s)$, respectively, and integrate from zero to infinity to obtain
\begin{align}
& u(0)\,a_1(\lambda)+\bar{u}_1(1+a_2(\lambda))+\bar{u}_2\,a_3(\lambda)=0,\label{u1eq}\\
& u(0)\,a_4(\lambda)+\bar{u}_1\,a_5(\lambda)+\bar{u}_2\,(1+a_6(\lambda))=0,\label{u2eq}
\end{align}
where
\begin{align*}
& a_1(\lambda)=-\int_0^{\infty} c(s)\,\alpha_1(s)\,\pi(s,\lambda)\,ds,\\
& a_2(\lambda)=\int_0^{\infty}c(s)\,\alpha_1(s)\,\pi(s,\lambda)\,\int_0^s\frac{g_1(r)}{\gamma^*(r)\,\pi(r,\lambda)}\,dr\,ds,\\
& a_3(\lambda)=\int_0^{\infty}c(s)\,\alpha_1(s)\,\pi(s,\lambda)\,\int_0^s\frac{g_2(r)}{\gamma^*(r)\,\pi(r,\lambda)}\,dr\,ds,\\
& a_4(\lambda)=-\int_0^{\infty} \alpha_2(s)\,\pi(s,\lambda)\,ds,\\
& a_5(\lambda)=\int_0^{\infty}\alpha_2(s)\,\pi(s,\lambda)\,\int_0^s\frac{g_1(r)}{\gamma^*(r)\,\pi(r,\lambda)}\,dr\,ds,\\
& a_6(\lambda)=\int_0^{\infty}\alpha_2(s)\,\pi(s,\lambda)\,\int_0^s\frac{g_2(r)}{\gamma^*(r)\,\pi(r,\lambda)}\,dr\,ds.
\end{align*}
Next we insert the solution \eqref{eigvsol} into the boundary condition \eqref{eigvbound} to obtain
\begin{equation}\label{boundeq}
u(0)\,(1+a_7(\lambda))+\bar{u}_1\,(a_8(\lambda)+g_3)+\bar{u}_2\,a_9(\lambda)=0,
\end{equation}
where
\begin{align*}
& a_7(\lambda)=-\int_0^{\infty}\frac{\beta(s)}{\gamma^*(0)}\,\pi(s,\lambda)\,ds,\\
& a_8(\lambda)=\int_0^{\infty}\frac{\beta(s)}{\gamma^*(0)}\,\pi(s,\lambda)\,\int_0^s\frac{g_1(r)}{\gamma^*(r)\,\pi(r,\lambda)}\,dr\,ds,\\
& a_9(\lambda)=\int_0^{\infty}\frac{\beta(s)}{\gamma^*(0)}\,\pi(s,\lambda)\,\int_0^s\frac{g_2(r)}{\gamma^*(r)\,\pi(r,\lambda)}\,dr\,ds.
\end{align*}
Now we can give a characterization of part of the point spectrum of the operator $\mathcal{A+B+C}$.
\begin{theorem}\label{eigenvalues}
For any $\lambda\in \left\{z\in\mathbb{C}\,|\, \text{\rm Re}\,z >-\mu_0\right\}$, we have
$\lambda\in\sigma_p(\mathcal{A+B+C})$ if and only if $\lambda$ satisfies the equation
\begin{equation}
K(\lambda)\eqdef \det 
\left( \begin{array}{llll}
a_1(\lambda) & 1+a_2(\lambda) & a_3(\lambda)\\
a_4(\lambda) & a_5(\lambda) & 1+a_6(\lambda)\\
1+a_7(\lambda) & g_3+a_8(\lambda) & a_9(\lambda)\\
\end{array}\right)=0.
\end{equation}
\end{theorem}
\Pf{.} If $\lambda\in\sigma_p(\mathcal{A+B+C})$, then Eqs.~\eqref{eigv}, \eqref{eigvbound} admit a nontrivial solution $u$. Hence for this $\lambda$ there exists a nonzero solution
vector $(u(0), \bar{u}_1, \bar{u}_2)$ such that Eqs.~\eqref{u1eq}-\eqref{boundeq} hold true. Thus $K(\lambda)=0$. Conversely, if $K(\lambda)=0$ for some $\lambda$ and $(u(0), \bar{u}_1, \bar{u}_2)$ is a nonzero solution of Eqs.~\eqref{u1eq}-\eqref{boundeq}, then $u$, given by Eq.~\eqref{eigvsol}, is a nonzero solution of Eqs.~\eqref{eigv}-\eqref{eigvbound} at least if $u(0)\not= 0$. If, however, $u(0)=0$, the only possible scenario for $u$ to vanish when defined by Eq.~\eqref{eigvsol} would be the condition that
\begin{equation}
	\bar{u}_1\,g_1\equiv -\bar{u}_2\,g_2
\end{equation}
hold true. Then Eqs.~\eqref{u1eq},\eqref{u2eq} would immediately give $\bar{u}_1=0=\bar{u}_2$ in contradiction to our assumption on $(u(0), \bar{u}_1, \bar{u}_2)$.
\eofproof

\noindent Since 
\begin{equation}
	\lim_{\lambda\rightarrow +\infty} K(\lambda)= \det 
\begin{pmatrix} 0 & 1 & 0\\
				0 & 0 & 1\\
				1 & g_3 & 0
\end{pmatrix} = 1,
\end{equation}
the limit being taken in $\mathbb R$, we can formulate the following simple instability criterion, which follows immediately from the Intermediate Value Theorem.

\begin{theorem}\label{instabeig}
A stationary solution $n_*$ of
Eqs.~\eqref{pde}--\eqref{density2} is linearly unstable if $K(0)<0$.
\end{theorem}
Note that for $n_*\equiv 0$
\begin{equation}
	K(0) =1-R(0).
\end{equation}
Hence the stationary solution $n_*\equiv 0$ is linearly unstable if $R(0)>1$.
In the remainder of this section we will concentrate on positive stationary solutions.

\example
Let us assume that 
the rate of an individual of size $s$ to attack another individual is proportional to the product of the probability for an individual of size $s$ to be attacked and its energetic value. Mathematically, this condition is modeled by
the relation
\begin{equation}
c(s)\alpha_1(s)=p\,\alpha_2(s),\quad s\in [0,\infty),\quad p\in\mathbb{R}^+. \label{alpha2}
\end{equation}
The constant $p$ denotes here the proportionality factor.
This condition is biologically relevant since environmental pressure conceivably makes individuals of higher energetic value, which are usually larger in size,  not only more prone to be attacked, but also more aggressive.
In addition we will assume that the attack rate of minimal size individuals (newborns) equals zero, i.e.\,$\alpha_2(0)=0$. With 
\begin{equation}
	\pi(s)\eqdef\pi(s,0)
\end{equation}
we obtain
\begin{equation}
	K(0)= \int_0^\infty \alpha_2(s)\,\pi(s)\,ds\,\int_0^\infty \frac{\beta(s)}{\gamma^*(0)}\,\pi(s)\,\int_0^s\frac{p\,g_1(r)+g_2(r)}{\gamma^*(r)\,\pi(r)}\,dr\,ds.
\end{equation}
Note that $\beta\not\equiv 0$ by Eq.~\eqref{netgrowth}.
Hence the instability criterion of Theorem~\ref{instabeig} is satisfied if, for $s\geq 0$,
\begin{align}
 & p^{-1}\,\alpha_1(s)+\alpha'_2(s)\,\gamma_E(s,E_*(s))+\alpha_2(s)\,
\left(
{d\over ds}\gamma_E(s,E_*(s))\right.\nonumber\\
&\left.-\gamma_E(s,E_*(s))\,
\left(\frac{\gamma^*_s(s)+\mu(s,E_*(s))+M_*(s)}{\gamma^*(s)}\right)
+\mu_E(s,E_*(s))\right)<0,\label{instabcond}
\end{align}
where we have used the relation
\begin{equation*}
n'_*(s)=-n_*(s)\,\left(\frac{\gamma^*_s(s)+\mu(s,E_*(s))+M_*(s)}{\gamma^*(s)}\right).
\end{equation*}
This instability condition automatically excludes the case $\alpha_2\equiv 0$.

\begin{remark}
In the age-structured case where $\gamma\equiv 1$ and where we may use $a$ for age in lieu
of size $s$ we find for the preceding example
\begin{equation}
M(a,t)=E(a,t)\,p^{-1}\,\frac{\alpha_1(a)}{\alpha_2(a)},
\end{equation}
assuming that $\alpha_2$ does not vanish. 
Hence the net reproduction function $R(n)$ -- as in the case of scramble competition -- can be considered a function of the environment $E$, i.e.
\begin{equation}\label{netrep2}
R(n)=\tilde{R}(E)=\int_0^{\infty}\beta(a)\,\exp\left\{-\int_0^a\mu(y,E(y,\cdot))+E(y,\cdot)\,
p^{-1}\,\frac{\alpha_1(y)}{\alpha_2(y)}\,dy\right\}\,da.
\end{equation}
When interpreting $\tilde{R}$ as a nonlinear operator between sets of bounded continuous functions  and ignoring
all issues pertaining to regularity, we may formally deduce the Fr\'echet derivative of $\tilde{R}$ at a stationary state $E_*$. The result of this formal calculation is
\begin{align} 
\tilde{R}_E(E_*)= &-\int_0^{\infty} \beta(a)\,\exp\left\{-\int_0^a\mu(y,E_*(y))+E_*(y)\,p^{-1}\frac{\alpha_1(y)}{\alpha_2(y)}\,dy\right\}\nonumber \\ 
&\quad \times\left(\int_0^a\mu_E(y,E_*(y))+p^{-1}\,\frac{\alpha_1(y)}{\alpha_2(y)}\,dy\right)\,da.\label{netrep2der}
\end{align}
Since in the scenario considered the instability condition \eqref{instabcond} reduces to 
\begin{equation}\label{specinstab}
	p^{-1}\,\alpha_1(a)+\alpha_2(a)\,\mu_E(a,E_*(a))<0,\quad a\in [0,\infty),
\end{equation}
Eq.~\eqref{netrep2der} implies
\begin{equation}
\tilde{R}_E(E_*)>0.
\end{equation} 
Hence in this special case of  model ingredients 
condition \eqref{instabcond} allows a formal, but intuitively clear biological interpretation, similar to the scramble competition case in \cite{FH1}: 
 if the net reproduction rate $\tilde{R}$ is increasing at a stationary environment $E_*$, then the equilibrium is unstable.
\end{remark}

\sector{Further spectral analysis}

Throughout this section we consider a positive stationary solution $n_*$ and assume that condition \eqref{muassume} holds true.

\begin{theorem}\label{specbound}
\begin{equation}\label{criterion}
\sigma(\mathcal{A+B})\cap \{z\in {\mathbb C}\,|\, \text{\rm Re}\,z> -\mu_0\}=P\cup \{0\},
\end{equation}
where the set $P\cup \{0\}$ consists of simple, isolated eigenvalues $\lambda$ of $\mathcal{A+B}$ such that
\begin{equation}
	-\mu_0< \text{\rm Re}\,\lambda<0\quad \text{for $\lambda\in P$.}
\end{equation}
\end{theorem}
\Pf{.}
Suppose first that $\lambda\in {\mathbb C}$ with $\text{Re}\,\lambda>  -\mu_0$ is such that
\begin{equation}
	\Lambda\left(\exp\left\{-\displaystyle{\int_0^\cdot\frac{\rho^*(r)+\lambda}{\gamma^*(r)}\,dr}\right\}\right)\not=1.
\end{equation}
Then the equation
\begin{equation}\label{res}
\left(\lambda\mathcal{I}-(\mathcal{A+B})\right)\,u = f
\end{equation}
with $f\in L^1(0,\infty)$ has the unique solution
\begin{align}
u(s)=&\ u(0)\,\exp\left\{-\int_0^s\frac{\rho^*(r)+\lambda}{\gamma^*(r)}\,dr\right\}\nonumber\\
& +\int_0^s\exp\left\{-\int_y^s\frac{\rho^*(r)+\lambda}{\gamma^*(r)}\,dr\right\}\frac{f(y)}{\gamma^*(y)}\,dy,\label{ressol}
\end{align}
where
\begin{equation}\label{resbound}
u(0)=\frac{\Lambda\left(\displaystyle{\int_0^\cdot\exp\left\{-\int_y^\cdot\frac{\rho^*(r)+\lambda}{\gamma^*(r)}\,dr\right\}\frac{f(y)}{\gamma^*(y)}\,dy}\right)}{1-\Lambda\left(\exp\left\{-\displaystyle{\int_0^\cdot\frac{\rho^*(r)+\lambda}{\gamma^*(r)}\,dr}\right\}\right)}.
\end{equation}
Note that our condition on $\lambda$ ensures that $\Lambda$ can be applied to deduce 
\eqref{resbound}. Moreover, it readily follows that $u$ belongs to $W^{1,1}(0,\infty)$. Hence we conclude that $\lambda\in \rho({\mathcal A+B})$. Now suppose that $\lambda\in {\mathbb C}$ with $\text{Re}\,\lambda> -\mu_0$ is such that
\begin{equation}\label{Lambdaeq}
	\Lambda\left(\exp\left\{-\displaystyle{\int_0^\cdot\frac{\rho^*(r)+\lambda}{\gamma^*(r)}\,dr}\right\}\right)=1.
\end{equation}
In light of Eq.~\eqref{netgrowth}, $\lambda=0$ solves \eqref{Lambdaeq}. In fact, $\lambda=0$
is the only solution with $\text{Re}\,\lambda\geq 0$. 
Let
\begin{equation}
	L(\lambda)\eqdef\Lambda\left(\exp\left\{-\displaystyle{\int_0^\cdot\frac{\rho^*(r)+\lambda}{\gamma^*(r)}\,dr}\right\}\right)-1
\end{equation}
for $\text{Re}\,\lambda>-\mu_0$. Then $L$ is analytic for $\text{Re}\,\lambda>-\mu_0$ and
\begin{equation}
L'(0)=-\int_0^{\infty}\frac{\beta(s)}{\gamma^*(0)}\,\exp\left\{-\int_0^s\frac{\rho^*(r)}{\gamma^*(r)}\,dr\right\}\,\int_0^s\frac{1}{\gamma^*(r)}\,dr\,ds<0.
\end{equation}
Hence $0$ is a pole of the resolvent operator of $\mathcal{A+B}$. Any other zero of $L$ gives rise to a pole of the resolvent operator of $\mathcal{A+B}$ since $L$ is analytic and nonconstant in the simply connected set $\text{Re}\,\lambda>-\mu_0$.
Finally, the representation of the resolvent operator, given through the solution $u$ in \eqref{ressol}, \eqref{resbound}, shows that the spectral projection associated with each pole of the resolvent operator is a rank one operator.
\eofproof

When the operator $\mathcal C$ is compact as in the special class of attack rates given in \eqref{newC}, we immediately obtain the following result.

\begin{corollary}\label{prop}
Suppose that the operator $\mathcal C$ is compact. Then 
$\sigma\left(\mathcal{A+B+C}\right)\cap \{z\in{\mathbb C}\,|\,\text{\rm Re}\,z>-\mu_0\}$
consists of isolated eigenvalues of $\mathcal{A+B+C}$ of finite multiplicity.
\end{corollary}

\sector{Asynchronous exponential growth}

The purpose of this section is to gain deeper insight into asymptotic properties of solutions of the linearized system \eqref{linear}-\eqref{linboundary}.
In particular, we are interested in solutions of the linearized problem which grow exponentially in time such that the proportion 
of individuals within any size range  compared to the total population approaches a limiting value as time tends to infinity, independently of the size distribution of the initial population. This phenomenon is called asynchronous exponential growth  and is known to be present, e.g.,  in the age-structured case.  Mathematical definitions will be given below.
The property of asynchronous exponential growth is important insofar as solutions can be regarded as asymptotically 
factorizable (with respect to time and size). Populations of this kind are often called ergodic (see \cite{I}). We refer to \cite{CH,NAG,GW,PL,WEB} for this and related notions.

In the framework of linear semigroup theory a strongly continuous semigroup $S=\left\{{\mathcal S}(t)\right\}_{t\geq 0}$ on a Banach space $\mathcal Y$
with generator ${\mathcal A}_S$ and growth bound
\begin{equation}
	s\left({\mathcal A}_S\right)\eqdef \sup\,\left\{\text{Re}\,\lambda\,|\, \lambda\in 
	\sigma\left({\mathcal A}_S\right)\right\}
\end{equation}
is said to exhibit balanced exponential growth (BEG for short) 
if there exists a bounded linear projection $\Pi$ on $\mathcal{Y}$ such that 
\begin{equation}\label{beg}
	\lim_{t\to \infty} \|e^{-s\left({\mathcal A}_S\right)\,t}\,\mathcal{S}(t)-\Pi\|=0.
\end{equation}
The semigroup $S=\left\{{\mathcal S}(t)\right\}_{t\geq 0}$ is said to exhibit asynchronous
exponential growth (AEG for short) if it exhibits BEG with a rank one projection $\Pi$.
For positive semigroups there exist well-known characterizations of BEG and AEG, see \cite{CH,NAG}. Our analytical approach will be guided toward these results.

\begin{theorem}\label{pos}
Suppose that, for every $y\geq 0$ and a.e.\,$s\geq 0$,
\begin{align}
& c(y)\,\alpha(y,s)\,\left((\gamma_E(s,E_*(s))\,n_*(s))_s+\mu_E(s,E_*(s))\,n_*(s)\right)
\nonumber\\
    &\quad +
c(y)\,D_2\alpha(y,s)\,\gamma_E(s,E_*(s))\,n_*(s)+\alpha(s,y)\,n_*(s)\leq 0,
\label{pos1} \\
&\beta(s)-\gamma_E(0,E_*(0))\,n_*(0)\,c(s)\,\alpha(s,0)\ge 0. \label{pos2}
\end{align}
Then the semigroup $\{{\mathcal T}(t)\}_{t\geq 0}$, generated by the operator ${\mathcal A} + {\mathcal B}+ {\mathcal C}$, is positive.
\end{theorem}
\Pf{.}
Condition \eqref{pos1} implies that the operator $\mathcal{C}$ is positive, hence we can restrict ourselves to the operator $\mathcal{A+B}$. Condition \eqref{pos2} ensures that the functional $\Lambda$ is nonnegative.  Consequently, the solution \eqref{ressol}, \eqref{resbound} of the resolvent equation \eqref{res} is well-defined  and nonnegative if $\lambda>0$ is sufficiently large.
\eofproof

\begin{corollary}\label{poscor}
Suppose that condition \eqref{pos2} holds true. Then the semigroup, generated by $\mathcal{A+B}$, is positive.
\end{corollary}

Let us recall a useful characterization of irreducibility on $L^1(\Omega,m)$ (see \cite{NAG}): A strongly continuous, positive semigroup $S=\left\{{\mathcal S}(t)\right\}_{t\geq 0}$ on the Banach lattice $\mathcal Y=L^1(\Omega,m)$
with generator ${\mathcal A}_S$
is irreducible if, for $f\in \mathcal Y$ with $f>0$, $\left(\lambda\,{\mathcal I}-{\mathcal A}_S\right)^{-1}\,f(s)>0$ for $m$-almost all $s\in \Omega$ and some $\lambda>s\left({\mathcal A}_S\right)$ sufficiently large.

\begin{theorem}\label{irr}
Suppose that the positivity conditions \eqref{pos1}, \eqref{pos2} hold true. Then the semigroup $T=\{\mathcal{T}(t)\}_{t\ge 0}$, generated by $\mathcal{A+B+C}$, is irreducible. 
\end{theorem} 
\Pf{.}
Since $C$ is positive and since, for $\lambda>0$  sufficiently large, the resolvent operator of $\mathcal A+B$ is positive as a consequence of Corollary~\ref{poscor}, we deduce from
\begin{equation}
	\left(\lambda\,{\mathcal I}-\left({\mathcal A+B+C}\right)\right)^{-1}=
	\sum_{n=0}^\infty \left(\left(\lambda\,{\mathcal I}-\left({\mathcal A+B}\right)\right)^{-1}\,{\mathcal C}\right)^n\,\left(\lambda\,{\mathcal I}-\left({\mathcal A+B}\right)\right)^{-1}
\end{equation} 
that it suffices to prove the irreducibility of the semigroup generated by $\mathcal A+B$. This result, however, follows immediately from the representation of solutions of the resolvent equation \eqref{res}, given by \eqref{ressol}, \eqref{resbound}.
\eofproof

Before we formulate the main result of this section, let us review the notions of essential norm, growth bound, and essential growth bound, and some of their properties.
Our discussion follows closely \cite{NAG}. 
Suppose that ${\mathcal A}_S$ is the infinitesimal generator of the strongly continuous semigroup $S=\left\{{\mathcal S}(t)\right\}_{t\geq 0}$ on a Banach space $\mathcal Y$. Then the growth bound of the semigroup is
defined by
\begin{equation}
	\omega_0\left({\mathcal A}_S\right)\eqdef \lim_{t\rightarrow \infty} {{\ln \|{\mathcal S}(t)\|}\over t}.
\end{equation}
For a bounded linear operator $T$ on $\mathcal Y$, the  essential norm
is given by
\begin{equation}
	\|T\|_\text{ess}\eqdef\text{dist}\,\left(T,{\mathcal K}({\mathcal Y})\right),
\end{equation}
where ${\mathcal K}({\mathcal Y})$ denotes the set of compact linear operators on $\mathcal Y$.
Of course, the essential norm is generally not a norm on the set of bounded linear operators on $\mathcal Y$. It is, however, a norm on the Calkin algebra of $\mathcal Y$, see \cite{NAG} and the references therein. Finally, the essential growth bound of the semigroup $S=\left({\mathcal S}(t)\right)_{t\geq 0}$ on $\mathcal Y$ with generator ${\mathcal A}_S$ is defined by
\begin{equation}\label{ess1}
	\omega_\text{ess}\left({\mathcal A}_S\right)\eqdef \lim_{t\rightarrow \infty} {{\ln \|{\mathcal S}(t)\|_\text{ess}}\over t}.
\end{equation}
It is readily seen that, for ${\mathcal K}\in {\mathcal K}({\mathcal Y})$,
\begin{equation}\label{ess2}
	\omega_\text{ess}\left({\mathcal A}_S\right)=\omega_\text{ess}\left({\mathcal A}_S+{\mathcal K}\right).
\end{equation}
The significance of the essential growth bound lies in the central fact that
\begin{equation}
	\omega_0\left({\mathcal A}_S\right)= \max\,\left\{\omega_\text{ess}\left({\mathcal A}_S\right), s\left({\mathcal A}_S\right)\right\}.
\end{equation}

\begin{theorem}\label{AEG}
Given a positive stationary solution $n_*$, suppose that  conditions \eqref{muassume}, \eqref{pos1}, and \eqref{pos2} hold true and that the operator $\mathcal C$ is  compact. If 
\begin{equation}\label{righthalf}
	\sigma\left(\mathcal{A+B+C}\right)\cap \{z\in{\mathbb C}\,|\,\text{\rm Re}\, z>0\}\not =\emptyset,
\end{equation}
then the linear semigroup $T=\{{\mathcal T}(t)\}_{t\geq 0}$, generated by
 $\mathcal{A+B+C}$, exhibits  AEG.
\end{theorem}
\Pf{.}
First we note that $\mathcal{A+B}$ has nonempty spectrum and generates a positive semigroup by Theorem~\ref{specbound} and  Corollary~\ref{poscor}. Hence Derndinger's Theorem (see \cite{De,NAG}) proves that
\begin{equation}
	\omega_0(\mathcal{A+B})=s(\mathcal{A+B})=0,
\end{equation}
where the last equality follows from Theorem~\ref{specbound}\,. Similarly, we obtain
from Theorem~\ref{pos} and Derndinger's Theorem that
\begin{equation}
	\omega_0(\mathcal{A+B+C})=s(\mathcal{A+B+C})>0.
\end{equation}
Here the last inequality is given by assumption \eqref{righthalf}.
Consequently, in light of Eqs.~\eqref{ess1}, \eqref{ess2}, we have
\begin{equation}
	\omega_\text{ess}(\mathcal{A+B+C})=\omega_\text{ess}(\mathcal{A+B})\leq
	\omega_0(\mathcal{A+B})=0<\omega_0(\mathcal{A+B+C}).
\end{equation}
Hence by Theorem~\ref{pos} and Theorem~\ref{irr}, the semigroup $T=\{{\mathcal T}(t)\}_{t\geq 0}$ is positive and irreducible with essential growth bound strictly smaller than its growth bound. The claim follows now immediately from Theorems 9.10 and 9.11 in 
\cite{CH}.
\eofproof

\sector{Conclusion}  In this work we have studied the asymptotic behavior of solutions of 
a linearized size-structured cannibalism model, recently introduced in \cite{D1}. 
The vital rates in this model depend  on a structuring variable (size), which takes values 
in an unbounded set, and on an infinite dimensional interaction variable (environment), describing the 
environmental feedback on individuals. Population models of this type are  notoriously difficult to analyze.
The reason for this difficulty is that the essential spectrum of the linearized operator is typically not empty, 
and even the point spectrum cannot be characterized in general via zeros of a characteristic function. 
The latter obstacle has already been observed in \cite{FH4} for a similar quasilinear hierarchically size-structured model. 
Therefore, analytical results, in particular with respect to the qualitative behavior of solutions, 
are rather rare in the literature, at least to our knowledge.
We would like to point out that the emphasis in the present work was to demonstrate how analytical techniques can be developed and used to treat qualitative questions of physiologically structured population models, where the structuring variable is unbounded and competition is incorporated through infinite dimensional interaction variables.

Here, using two different strategies, we have formulated linear stability and instability
criteria for equilibrium solutions of the model.
We derived a general instability criterion in the case of a separable attack rate and extended our dissipativity
approach, employed in \cite{FH4} for a model with finite size span, to the case of infinite size span.
We also carried out a more refined spectral analysis of the linearized operator 
which allowed us to gain deeper insights into the asymptotic behavior of solutions of the linearized system.
In particular, we investigated the question whether solutions of the linearized problem
exhibit asynchronous exponential growth and gave sufficient conditions for an affirmative answer. There are only few results of this type 
for models with unbounded structuring variable, see \cite{PL,PR1,TH1,TH2}.
In passing, we also note that our spectral analysis gives  rise to an extension of our stability results in \cite{FH1} for a  model with finite size span to one with infinite size span.

\bigskip

{\bf Acknowledgment}

JZF was supported by the EPSRC grant EP/F025599/1. TH gratefully acknowledges
support through NSF-Grant DMS 0709197. We thank the reviewers for their thoughtful comments and suggestions.

\end{document}